\documentclass[preprint]{elsarticle}

\usepackage{amssymb}
\usepackage{amsfonts}
\usepackage{amsmath}
\usepackage{graphicx}
\usepackage{amsthm}
\usepackage{algorithm}
\usepackage{epic}

\newtheorem{rem}{Remark}

\newtheorem{exmpl}{Example}

\def\diag{\mathop{\mathrm{diag}}}
\def\trace{\mathop{\mathrm{trace}}}

\begin{document}

\begin{frontmatter}
 \title{Fast computation of optimal damping parameters for linear vibrational systems}

\author[FESB]{N.\ Jakov\v{c}evi\'{c} Stor\corref{c1}\fnref{f1}}
\ead{nevena@fesb.hr}
\author[FESB]{I.\ Slapni\v{c}ar\fnref{f1}}
\ead{ivan.slapnicar@fesb.hr}
\author[MATHOS]{Z.\ Tomljanovi\'{c}\fnref{f1}}
\ead{ztomljan@mathos.hr}

\cortext[c1]{Corresponding author}
\fntext[f1]{This work was supported by the Croatian Science Foundation  projects `Optimization of Parameter Dependent Mechanical Systems' (IP-2014-09-9540) and `Vibration Reduction in Mechanical Systems' (IP-2019-04-6774).}

\address[FESB]{Faculty of Electrical Engineering, Mechanical Engineering and
Naval Architecture, University of Split, Rudjera Bo\v{s}kovi\'{c}a
32, 21000 Split, Croatia}

\address[MATHOS]{Department of Mathematics, J. J. Strossmayer University of Osijek,

Trg Ljudevita Gaja 6, 31000 Osijek, Croatia}

\begin{abstract}
We formulate the quadratic eigenvalue problem underlying the mathematical model of a linear
vibrational system as an eigenvalue problem of a diagonal-plus-low-rank matrix $A$. The
eigenvector matrix of $A$ has a Cauchy-like structure. Optimal viscosities are those for which $\trace(X)$ is minimal, where $X$ is the solution of the  Lyapunov equation $AX+XA^{*}=GG^{*}$.  Here $G$ is a low-rank matrix which depends
on the eigenfrequencies that need to be damped.
After initial eigenvalue decomposition of linearized problem which requires $O(n^3)$ operations, our algorithm computes optimal viscosities for each choice of external dampers
in $O(n^2)$ operations, provided that the number of dampers is small.
Hence, the subsequent optimization is order of magnitude faster than in the standard approach which solves Lyapunov equation in each step, thus requiring $O(n^3)$ operations. Our algorithm is based on $O(n^2)$ eigensolver for  complex symmetric diagonal-plus-rank-one matrices and fast $O(n^2)$ multiplication of linked Cauchy-like matrices.
\end{abstract}


\end{frontmatter}

\section{Introduction}

We consider determination of optimal damping for the linear vibrational
system described by

\begin{equation}  \label{MDK}
M\ddot{x}+D\dot{x}+Kx=0,
\end{equation}
where $M$ and $K$ (called mass and stiffness, respectively) are real,
symmetric positive definite matrices of order $n$.
The damping matrix is defined as
\begin{equation} \label{D}
D=D_{int}+D_{ext},
\end{equation}
where matrices $D_{int}$ and $D_{ext}$ correspond to internal and external
damping, respectively.
Internal damping can be modeled as
\begin{equation} \label{Dint1}
D_{int}=\alpha \, M^{1/2}\sqrt{M^{-1/2}KM^{-1/2}}M^{1/2}, \quad \alpha\in
(0.005,0.1),
\end{equation}
or
\begin{equation}\label{Dint2}
D_{int}=\alpha M+\beta K,\quad \alpha ,\beta \in (0.005,0.1) .
\end{equation}

From the optimization point of view we are more interested in the external
viscous damping which can be modeled as
\begin{equation}
D_{ext}=\rho_{1}D_{1}+\rho_{2}D_{2}+\cdots +\rho_{k}D_{k},  \label{Dext}
\end{equation}
where $\rho_{i}$ is the viscosity and $D_{i}$ describes a geometry of the
corresponding damping position for $i=1,\ldots ,k$. Typically, system has
a very few dampers compared to the full dimension, which means that
$k\ll n$. More details on  the structure of internal and external damping can be found in e.g. \thinspace \cite{BTT13, BTT11, KTT12, Ves90, Ves11}.
The model of linear vibrational system \eqref{MDK} corresponds to the quadratic eigenvalue problem
\begin{equation}\label{QEP}
(\lambda^2 M + \lambda D + K)\,x=0,\quad   x\neq 0.
\end{equation}

The damping optimization problem, in general, can be stated as follows: determine the "best" damping matrix
$D$ which insures optimal evanescence of each component of $x$. There exists several optimization criteria for
this problem. One criterion is the so-called spectral abscissa criterion. This criterion requires that the maximal
real part of the
eigenvalues of the quadratic eigenvalue problem \eqref{QEP} are minimized
(see  e.g.\ \cite{FreLan99, NakTomTr13}).

In this paper we will use criterion  based on total average energy of the considered system.
This criterion considers minimization of the total energy of the system (as a sum of kinetic and potential energy)
averaged
over all initial states of the unit total energy and a given frequency range.
Benefits of this criterion are discussed in \cite{CoxNakRittVes04, Ves11}. Moreover, it can be shown
(see e.g.\ \cite{CoxNakRittVes04, Nak02, Ves11}) that this
minimization criterion is equivalent to
\begin{equation*}
\mathop{\mathrm{trace}}{X}\quad \rightarrow \quad \min ,
\end{equation*}
where $X$ is the solution of the Lyapunov equation
\begin{equation}
AX+XA^{\ast }=-GG^{T} \label{Lyap}
\end{equation}
where matrix $A$ comes form linearization of the system \eqref{MDK} and matrix $G$ determines
which part of an undamped eigenfrequencies needs to be damped.

In practice, one can usually influence only the external damping.
The problem is to determine the best damping matrix $D$, 
which minimizes the total
average energy of the system, that is, to determine the optimal matrix $D_{ext}$
for which $\mathop{\mathrm{trace}}X$ is minimal with $A$ as in \eqref{A}.

The paper is organized as follows: symmetric linearization and
an overview of the existing $O(n^{3})$ solution method, which minimizes trace of the solution of the respective Lyapunov equation,
are presented in Section 2. A new $O(n^2)$ algorithm, which reduces the problem to sequence of eigenproblems of complex symmetric diagonal-plus-rank-one (CSymDPR1) matrices and uses fast multiplication of linked Cauchy-like  matrices, is presented in Section 3. Numerical examples and some timings are given in Section 4.

\section{Standard solution}

In this section, we describe the standard solution method for finding the optimal damping of the linear vibrational system \eqref{MDK}. In Section 2.1,  we describe a linearization of the  the quadratic eigenvalue problem \eqref{QEP}, which is, by changing basis, further reduced to an eigenvalue problem of a simpler matrix.
In Section 2.2, we present
an existing direct approach (see e.g.\  \cite{CoxNakRittVes04, Ves11}), which is
$O(n^{3})$ solution of the problem.

For structured system this problem was considered in \cite{BTT13, BTT11} where authors
proposed dimension reduction to accelerate the optimization process. However, to be efficient, dimension reduction requires specific structure, so this approach cannot be applied efficiently in general setting.

\subsection{Symmetric linearization}\label{symlin}

By symmetric linearization we transform quadratic eigenvalue
problem \eqref{QEP}
to the
generalized eigenvalue problem (GEVP)  (see e.g.  \cite{TiMe01})
\begin{equation}\label{GEVP}
\begin{bmatrix}
0 & K \\
K & D%
\end{bmatrix}%
\begin{bmatrix}
y \\
x%
\end{bmatrix}
=\lambda%
\begin{bmatrix}
K &  \\
& -M%
\end{bmatrix}%
\begin{bmatrix}
y \\
x%
\end{bmatrix}%
.
\end{equation}

Let
\begin{equation*}
\Phi^T K \Phi = \Omega^2, \quad \Phi^T M \Phi =I .
\end{equation*}
the generalized eigenvalue decomposition of the pair $(K,M)$.
Since the calculation of matrices  $\Phi$ and $\Omega$ does not depend
on the damping matrix $D$, they can be  calculated
prior to the optimization procedure.

Both choices of $D_{int}$ from \eqref{Dint1} or \eqref{Dint2} imply that $D_{int}$ is diagonal in the
$\Phi$-basis, that is
\begin{equation*}
\Phi^T D_{int} \Phi=\Gamma=\mathrm{diag}.
\end{equation*}
The external damping matrix, $D_{ext}$, is a low-rank matrix of rank $k$ which depends on the number, positions and the structure of dampers.
For example, if all dampers are grounded and $l$ is a vector of indices of
damping positions, then $D_{ext}$ is zero except for $[D_{ext}]_{l[i],l[i]}=1$.

In the basis
$$\begin{bmatrix}
\Phi &  \\
& \Phi%
\end{bmatrix},$$
problem \eqref{GEVP} reduces to GEVP
\begin{equation*}
\begin{bmatrix}
0 & \Omega^2 \\
\Omega^2 & \Gamma + \Phi^T D_{ext} \Phi
\end{bmatrix} x=\lambda	
\begin{bmatrix}
\Omega^2 & 0 \\
0 & -I%
\end{bmatrix} x%
,
\end{equation*}
\bigskip
and in the basis
$$\begin{bmatrix}
\Phi\Omega^{-1} &  \\
& \Phi%
\end{bmatrix}%
$$
we have the hyperbolic generalized eigenvalue problem
\begin{equation}\label{HEVP}
\begin{bmatrix}
0 & \Omega \\
\Omega & \Gamma + \Phi^T D_{ext} \Phi%
\end{bmatrix}
x=\lambda	
\begin{bmatrix}
I & 0 \\
0 & -I%
\end{bmatrix}x.
\end{equation}

Now we can write the linearized system in the so-called modal coordinates.
By simple transformation, \eqref{HEVP} is equivalent to the eigenvalue problem for the matrix
\begin{equation}
A=
\begin{bmatrix} \label{A}
0 & \Omega  \\
-\Omega  & -(\Gamma +\Phi ^{T}D_{ext} \Phi )%
\end{bmatrix}.
\end{equation}
This is also the matrix $A$ from the Lyapunov equation \eqref{Lyap}.
The matrix $G$ from \eqref{Lyap} is equal to
\begin{equation} \label{G}
G=%
\begin{bmatrix}
I_s & 0 \\
0 & 0 \\
0 & I_s \\
0 & 0%
\end{bmatrix},
\end{equation}
where parameter $s$ determines the number of eigenfrequencies of the undamped system which have to be damped
(for more details see e.g.\ \cite{BTT13, TV09, Ves11}).

\subsection{Existing approach, one of, $O(n^{3})$ solution}

This is demanding problem,  both, from the computational point of view and from the point of optimization of  damping positions. The main reason lies in the fact that the criterion of total average energy has
many local minima, so we usually need to optimize viscosity parameters for many different damping positions.

Lyapunov equation (\ref{Lyap}) with structured matrices $A$ and $G$ from (\ref{A}) and (\ref{G}), respectively,
can be solved by iterative approaches such as ADI method
\cite{LW02, Penzl} used in \cite{TV09} or the sign function method \cite{BQ99, KL95} used in \cite{BD16}.

Standard direct approaches calculate the solution of Lyapunov equation \eqref{Lyap} by using  the Schur form, for example,  Hammarling algorithm \cite{Hammarling:1982:NSS, Kressner08} and Bartels-Stewart
algorithm \cite{BARTSTEW72}. The computation of Schur form requires $O(n^{3})$ operations, so these algorithms are
$O(n^{3})$ solutions. The algorithms are implemented in the \texttt{SLICOT} library
\cite {Slicot} and are used in \texttt{Matlab}.
The timings for some examples are given in Section 4.

\section{The new algorithm, $O(n^{2})$ solution}

In our approach, instead of using the Schur form, we use diagonalization of the matrix $A$ from (\ref{A}).
The eigenvalue problem for the matrix $A$ is reduced to $k$ eigenvalue problems for the
complex symmetric diagonal-plus-rank-one (CSymDPR1) matrices, $k$ being the number of dampers. Each of those EVPs can be efficiently solved in $O(n^2)$ operations.
It is important that updating of the eigenvectors can also be performed using $O(n^2)$ operations, due
to Cauchy-like structure of eigenvector matrices.
In this way, after preparatory steps from Sections \ref{symlin} and \ref{reduction} below, which require $O(n^3)$ operations,
each  computation of $\trace(x)$, where $X$ is from (\ref{Lyap}), requires only $O(n^2)$ operations. This makes trace optimization considerably faster.

The section is organized as follows. In Section \ref{cauchylike}, we present existing results about Cauchy-like
matrices and their fast multiplication.
In Section \ref{CSymDPR1} we develop an efficient $O(n^2)$ method for the solution of the CSymDPR1 eigenvalue problem. In Section \ref{reduction}, we describe reduction to the CSymDPR1 eigenvalue problems.
In Section \ref{Trace}, we develop fast $O(n^2)$ algorithm for the final trace computation, based on the fast multiplication of Cauchy-like matrices.

\subsection{Cauchy-like matrices}\label{cauchylike}

Cauchy-like matrix $C(X,Y,P,Q)$ is the matrix which satisfies the Sylvester-type displacement equation
(see e.g.\ \cite{PaZh00})
\begin{equation}
XC-CY=P\cdot Q^{T}, \label{Cauchy-like}
\end{equation}
where $X=\mathop{\mathrm{diag}}(x_{l})$, $Y=\mathop{\mathrm{diag}}(y_{l})$, $l=1,\ldots, n$, and $P,Q \in \mathbb{R}^{n\times k}$. The matrices $X$, $Y$, $P$ and $Q$ are called the generators of $C$.
For example, the standard Cauchy matrix $C=[1/(x_i-y_j]$ is equal to
$C(X,Y,\mathbf{1},\mathbf{1})$, where $\mathbf{1}=\begin{bmatrix}
1 & \ldots & 1 \end{bmatrix}^{T}$. Clearly, given generators, all elements of a Cauchy-like matrix can be computed on $O(kn^2)$ operations.

For multiplication by Cauchy-like matrices, we have following results.

Given $n\times n$ Cauchy-like matrix $A$ from \eqref{Cauchy-like} and $n$-dimensional vector $v$, the product $Av$ can be computed
 in $O(nk\log n^{2})$ operations \cite[Lemma 2.2]{PaZh00}.

Given two linked Cauchy-like matrices such that
\begin{equation}\label{clm}
\mathop{\mathrm{diag}}(a)\, A - A \mathop{\mathrm{diag}}(b)=EF^*, \quad %
\mathop{\mathrm{diag}}(b)\, B - B \mathop{\mathrm{diag}}(c)=NH^*,
\end{equation}
where $E,F \in \mathbb{R}^{n\times k_{1}}$ and $N,H \in \mathbb{R}^{n\times k_{2}}$,
the product $C=A\cdot B$ is a Cauchy-like matrix which satisfies \cite[Lemma 2.3]{PaZh00}

\begin{equation}\label{clp}
\mathop{\mathrm{diag}}(a)\, C - C \mathop{\mathrm{diag}}(c)=
\begin{bmatrix}
E & AN%
\end{bmatrix}
\begin{bmatrix}
B^*F & H%
\end{bmatrix}%
^*.
\end{equation}
This generators of $C$ can be computed in $O(nk_{1}k_{2} \log n^{2})$ operations.

\subsection{Eigenvalue decomposition of CSymDPR1 matrix}\label{CSymDPR1}

Let $A$ be an $n\times n$ CSymDPR1 matrix,
\begin{equation}
A= \Xi +\rho y y^T,
\label{Adpr1}
\end{equation}%
where $\Xi=\mathop{\mathrm{diag}}(\xi_{1},\xi_{2},\ldots ,\xi_{n})\text{ }$
is a diagonal matrix,
$y=\begin{bmatrix}y_{1} & y_{2} &\cdots  &y_{n}\end{bmatrix}^{T}$ is a vector,
 and $\rho \neq 0$ is a real scalar.
Here $\xi_i,y_i\in\mathbb{C}$.

Without loss of generality, we assume that
$A$ is irreducible, that is,
\begin{align*}
y _{i}&\neq 0,\quad i=1,\ldots, n, \\
\xi_{i}&\neq \xi_{j},\quad i\neq j,\quad i,j=1,\ldots ,n.
\end{align*}
Indeed, if $y_{i}=0$ for some $i$, then the diagonal element $\xi_{i}$ is an
eigenvalue whose corresponding eigenvector is the $i$-th canonical vector, and
if $\xi_{i}=\xi_{j}$, then $\xi_{i}$ is an
eigenvalue of the matrix $A$.

Let
\begin{equation*}
A=V\Lambda V^{T}  \label{Aeigendec}
\end{equation*}%
be the eigenvalue decomposition of $A$, where $\Lambda =\mathop{\mathrm{diag}}(\lambda _{1},\lambda _{2},\ldots ,\lambda
_{n})$ are the eigenvalues and
$V=\left[
\begin{array}{ccc}
v_{1} & \cdots & v_{n}%
\end{array}%
\right]
$ is a matrix whose columns are the corresponding eigenvectors. Notice that the eigenvector matrix of  a complex symmetric matrix satisfies the relation $V^{-1}=V^T$.

The eigenvalue problem for $A$ can be solved by any of the
standard methods (see \cite{Wat14} and the references therein).
However, due to the special structure of $A$, we
can use the following approach (see \cite{Cup81}  and \cite[Section 8.5.3]{GV96}):
the eigenvalues of $A$ are the zeros of the secular equation
\begin{equation}
1+\rho\sum_{i=1}^{n}\frac{y_{i}^{2}}{
\xi_{i}-\lambda }\equiv 1 +\rho y^{T}(\Xi-\lambda I)^{-1}y=0,  \label{Pick}
\end{equation}
and the corresponding eigenvectors are given by
\begin{equation}
v_{i}=\frac{x_{i}}{\left\Vert x_{i}\right\Vert _{2}},\text{ \ \ \ }
x_{i}=( \Xi-\lambda _{i}I) ^{-1}y,\quad i=1,\ldots ,n.  \label{Aeigenvec}
\end{equation}
It is important to notice that $V$ is a Cauchy-like matrix,
\begin{equation}\label{VCL}
V=C(\Xi,\Lambda,y,\Psi),\quad  \Psi=\begin{bmatrix}\displaystyle\frac{1}{\|x_1\|_2} & \cdots &\displaystyle\frac{1}{\|x_n\|_2}\end{bmatrix}^T.
\end{equation}
The equation \eqref{Pick} can, for example, be solved by the secular equation solver from the package MPSolve package  \cite{BiFi00, MPSolve}.

If $A$ is real, the eigenvalues interlace the diagonal elements of $\Xi$, and
can be computed highly accurately by bisection \cite{JSB13}. In this case,
orthogonality of computed eigenvectors follows from the accuracy of computed $\lambda$s.
In the complex symmetric case there is no interlacing, but orthogonality is
not an issue, so we developed a version of the Rayleigh quotient iteration.

Standard Rayleigh quotient iteration (RQI) is as follows \cite{Par74}: given starting $x$ repeat

\begin{equation}\label{RQI}
\mu=\frac{x^* Ax}{x^*x}, \quad x:=(A-\mu I)^{-1}x.
\end{equation}
Than $\mu\to \lambda$. In our case
\begin{equation*}
(A-\mu I)^{-1}=\Xi^{-1}+\gamma (\Xi-\mu I)^{-1}yy^{T}(\Xi-\mu I)^{-1},\quad \gamma =-\frac{\rho }{1+\rho
y^{T}(\Xi-\mu I)^{-1}y}.
\end{equation*}
is again a CSymDPR1 matrix which is computed in $O(n)$ operations.

For real symmetric or Hermitian matrices, RQI converges quadratically to absolutely largest eigenvalue. In the complex symmetric case, the convergence of RQI is slow and it is better to use the modified Rayleigh quotient iteration (MRQI) which is as follows: given starting $x$ repeat
\begin{equation}
\mu=\frac{x^T Ax}{x^T x}, \quad x:=(A-\lambda I)^{-1}y. \label{MRQI}
\end{equation}
The MRQI method converges quadratically \cite{ArHo04,Par74}.

For a CSymDPR1 matrix, having in mind the eigenvector formulas \eqref{Aeigenvec}, we further modify the method as follows: given starting $x$ repeat
\begin{equation}
\mu=\frac{x^T Ax}{x^T x}, \quad x:=(\Xi-\lambda I)^{-1}y. \label{MRQI1}
\end{equation}
This modified method showed very good convergence properties in all our large damping problems.

Once $\mu $ has converged to an eigenvalue, this eigenvalue can be deflated \cite[Section 7.2]{PaZh11}.
In particular,  if for some $l<n$ we have computed eigenvalues $\lambda_{n-l+1},\ldots,\lambda_{n}$ of $A$, then we can compute the remaining $n-l$ eigenvalues
$\lambda_{1},\ldots,\lambda_{n-l}$ as eigenvalues of the $(n-l) \times (n-l)$ CSymDPR1 matrix
$$
\widetilde A = \widetilde{\Xi} +\rho \widetilde{y}\widetilde{y}^T,
$$
where $\widetilde{\Xi}=\diag (\xi_{1},\xi_{2},\ldots,\xi_{n-l})$ and
\begin{equation*}
\widetilde{y}_{i}=y_{i}\prod\limits_{j=1}^{n-l}\frac{\xi_{i}-\xi_{n+1-j}}{%
\xi_{i}-\lambda _{n+1-j}},\quad i=1,\ldots,n-l.
\end{equation*}

In our implementation, first steps use RQI from \eqref{RQI} and, after that, MRQI from \eqref{MRQI1} is used until convergence.

The operation count to compute all eigenvalues of $A$ is $O(n^2)$, construction of generators for the eigenvector matrix $V$ from \eqref{VCL}  takes $O(n^2)$ operations (computing $\Psi$), and the reconstruction of $V$ from its generators, if needed, takes another $O(n^2)$ operations. This amounts to $O(n^2)$ operations to compute the complete eigenvalue decomposition of $A$.

\subsection{Reduction to CSymDPR1 eigenproblems}\label{reduction}

Let $\Xi$ and $Q$ denote the solution of the hyperbolic GEVP

\begin{equation}\label{HEVP1}
\begin{bmatrix}
0 & \Omega \\
\Omega & \Gamma
\end{bmatrix}x
=\lambda	
\begin{bmatrix}
I & 0 \\
0 & -I
\end{bmatrix}
,
\end{equation}
where $\Omega$ and $\Gamma$ are defined by \eqref{HEVP}, such that
\[
Q^T\begin{bmatrix}
0 & \Omega \\
\Omega & \Gamma
\end{bmatrix} Q =\Xi, \quad Q^T\begin{bmatrix}
I & 0 \\
0 & -I
\end{bmatrix}Q=I.
\]
Due to sparse structure of the problem \eqref{HEVP1}, the matrices $Q$ and $\Xi$ are computed by solving $n$ eigenvalue problems for $2\times 2$ matrices: for
$i=1,2,\ldots,n$,
$$
\begin{bmatrix}Q_{ii} & Q_{i,i+n} \\ Q_{i+n,i} & Q_{i+n,i+n}\end{bmatrix}^{-1}
 \begin{bmatrix} 0 & \Omega_{ii} \\ -\Omega_{ii} & \Gamma_{ii}\end{bmatrix} \begin{bmatrix}Q_{ii} & Q_{i,i+n} \\ Q_{i+n,i} & Q_{i+n,i+n}\end{bmatrix}= \begin{bmatrix} \Xi_{ii} & 0 \\ 0 & \Xi_{n+i,n+i}\end{bmatrix},
$$
and all other elements of $Q$ and $\Xi$ are zero.

The problem \eqref{HEVP1} is equal to the problem \eqref{HEVP}, but without external damping. Instead of solving \eqref{HEVP}, we now compute the eigenvalue decomposition of the complex symmetric diagonal-plus-low-rank matrix
\begin{equation}\label{dprk}
A= \Xi+ Q^T \begin{bmatrix} 0 & 0 \\ 0 & \Phi^T D_{ext}\Phi \end{bmatrix}Q.
\end{equation}
Assume, for example, that there is only one damper with viscosity $\rho$ positioned at the mass $l$.
Instead of solving \eqref{HEVP}, we compute the eigenvalue decomposition of the CSymDPR1 matrix
\begin{equation*}
A= \Xi+ \rho yy^T,
\end{equation*}
where $y=Q_{n+1:2n,:}^T\Phi_{l,:}$.
In the case of $k$ dampers, the procedure is repeated.
For example, in the case of two dampers we need to solve the eigenproblem for the matrix
$$A=\Xi +\rho _{1}y_{1}y_{1}^{T}+\rho _{2}y_{2}y_{2}^{T}.$$
First, we find eigendecomposition of matrix
$$\Xi +\rho _{1}y_{1}y_{1}^{T}=S_{1}\Lambda_1 S_{1}^{T}.$$
Then the eigendecomposition of $A$ is computed as
\begin{align*}
A&=S_{1}\Lambda_1 S_{1}^{T}+\rho _{2}y_{2}y_{2}^{T} \\
&=S_{1}(\Lambda_1 +\rho _{2}S_{1}^{-1}y_{2}y_{2}^{T}S_{1}^{-T})S_{1}^{T}\\
&=S_1 S_2 \Lambda	S_2^T S_1^T\equiv S\Lambda S^T.
\end{align*}
Since $A$ is complex symmetric, we also have
\begin{equation}\label{orth}
S_1^{-1}=S_1^T,\quad S_2^{-1}=S_2^T,\quad S^{-1}=S^T.
\end{equation}
From \eqref{VCL} it follows that $S_1$ and $S_2$ are Cauchy-like matrices,
$$
S_1=C(\Xi, \Lambda_1,y,\Psi_1),\quad S_2=C(\Lambda_1,\Lambda,S_1^T y_2,\Psi_2),
$$
where elements of vectors $\Psi_1$ and $\Psi_2$ are reciprocals of the norms of unnornmalized eigenvectors $x_i$ from \eqref{Aeigenvec}. The matrices $S_1$ and $S_2$ are linked, so according to \eqref{clm} and \eqref{clp}, $S=S_1\cdot S_2$ is a Cauchy-like matrix,
$$
S=C(\Xi,\Lambda,P,Q), \quad P=\begin{bmatrix} y &S_1 y_2\end{bmatrix},\quad
Q=\begin{bmatrix} S_2^T \Psi_1 &  \Psi_2\end{bmatrix}.
$$
This procedure is easily  generalized to $k>2$ dampers.

For small $k$, the computation of $\Lambda$, $P$ and $Q$ requires $O(n^2)$ operations.

\subsection{Trace computation} \label{Trace}

Let $A$ be given by \eqref{dprk} and let $A=S\Lambda S^T$ be its eigenvalue decomposition computed with the method from Section \ref{CSymDPR1}. Then
$S$ is a Cauchy-like matrix defined by $S=C(\Xi,\Lambda,P,Q)$ for some
$P,Q\in\mathbb{C}^{n\times k}$ satisfying $S^{-1}=S^T$, where $k$ is the number of dampers.

Let $\bar A$ denote the elementwise conjugated matrix $A$. Inserting the eigenvalue decomposition of $A$ into the Lyapunov equation  \eqref{Lyap} gives
$$
S\Lambda S^T X + X(S\Lambda S^T)^*=S\Lambda S^T X + X\bar S\bar \Lambda S^*=-GG^T $$
Premultiplying this equation by $S^T=S^{-1}$, postmultiplying by $\bar S=S^{-*}$ and setting $Y=S^TX\bar S$, gives a displacement equation
\begin{equation}
\Lambda  Y + Y \bar \Lambda
=-S^{T}GG^T  \bar S. \label{DisplY}
\end{equation}
Therefore, $Y$ is a Cauchy-like matrix, $Y=C(\Lambda, \bar\Lambda, -S^TG,S^T G)$.
Notice that $S^T G$ is not an actual matrix multiplication -- due to the special form of $G$ from \eqref{G}, this is just a selection of columns of $S^T$. Generating full $Y$, if needed, requires $O(sn^2)$ operations.

To finish the computation, we need to compute $\mathop{\mathrm{trace}}(X)$.
Set $Z=SY$. Then
$$\mathop{\mathrm{trace}}(X)=\mathop{\mathrm{trace}}(SYS^*)=
\mathop{\mathrm{trace}}(ZS^*)=\mathop{\mathrm{trace}}(S^*Z).
$$
Since $S$ and $Y$ are linked Cauchy-like matrices, according to
\eqref{clp}, $Z$ is a Cauchy-like matrix
$$
Z=S\cdot Y=C(\Xi,\Lambda,P,Q)\cdot C(\Lambda,\bar\Lambda,-S^TG,S^TG)=
C(\Xi, \bar\Lambda,P',Q'),
$$
where
$$
P'=\begin{bmatrix} P & -S S^T G\end{bmatrix}=
\begin{bmatrix} P & -G\end{bmatrix},\quad
Q'= \begin{bmatrix}Y^*Q & S^T G \end{bmatrix}.
$$
 The computation of $Q'$ requires $O(nks\log^2 n)$ operations (see Section \ref{CSymDPR1}).

  Finally, $\mathop{\mathrm{trace}}(X)$ is computed using scalar products:
 $$
  \mathop{\mathrm{trace}}(X)=  \mathop{\mathrm{trace}}(S^* Z)=
  \sum_{i=1}^n  \prod_{k=1}^n  \bar S_{ik}Z_{ki}.
   $$
   This step requires $O(kn^2)$ operations to compute elements of $\bar S$,
   $O((k+s)n^2)$ operations to compute elements of $Z$, and $O(n^2)$ operations to compute scalar products.


\section{Numerical examples}

In this section we present three examples of $n$-mass oscillator.
The size of the problem is $n=801$ for the ``small'' example,  $n=1601$ for the ``large'' example, and
$n=2001$ for the ``homogeneous'' example with more homogeneous masses.
We compare our algorithm with the $O(n^{3})$ algorithm from \cite{BTT13}.
The computations are performed on an Intel i7-8700K CPU running at 3.70GHz with 12 cores.

\begin{figure}[ht]
\begin{center}\centering
\includegraphics[width=9cm]{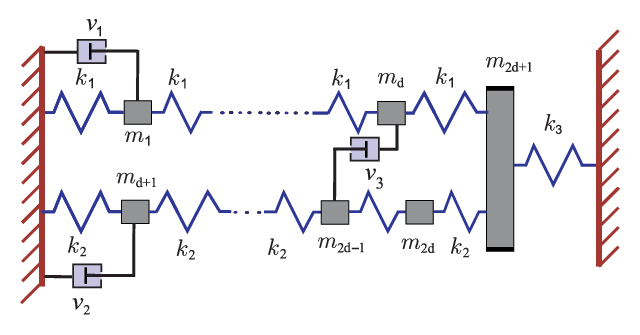}
\caption{$2d+1$ mass oscillator}
\label{2doscillator}
\end{center}
\end{figure}

Let us describe the large example. The small example is similar. We consider the mechanical system shown in
Figure~\ref{2doscillator}. The mass oscillator contains two rows of $d$ masses that
are grounded from one side, while on the other side masses are connected to
one mass which is then grounded. Therefore, we consider $2d+1$ masses and
$2d+3$ springs, while the system has three dampers of different viscosities
$\rho_1$, $\rho_2$ and $\rho_3$.

The mass matrix is
\begin{equation} \label{M}
M=\mathop{\mathrm{diag}}
(m_1,m_2,\ldots,m_n),
\end{equation}
and the stiffness matrix is
\begin{equation} \label{K}
K=%
\begin{bmatrix}
K_{11} &  & -\kappa _{1} \\
& K_{22} & -\kappa _{2} \\
-\kappa _{1}^{T} & -\kappa _{2}^{T} &  k_{1}+k_{2}+k_{3}%
\end{bmatrix}
,
\end{equation}
where
\begin{equation*}
K_{ii}=k_{i}
\begin{bmatrix}
2 & -1 &  &  &  \\
-1 & 2 & -1 &  &  \\
& \ddots & \ddots & \ddots &  \\
&  & -1 & 2 & -1 \\
&  &  & -1 & 2
\end{bmatrix}
,\quad \kappa _{i}=
\begin{bmatrix}
0 \\
\vdots \\
0 \\
k_{i}
\end{bmatrix}, \quad i=1,2
.
\end{equation*}
For $d=800$ and $n=1601$ masses, we consider the
following configuration:
\begin{align*}
k_{1}& =100,\quad k_{2}=150,\quad k_{3}=200,  \\
m_{i}& =2000-4i, \quad i=1,\ldots ,d/2, \\
m_{i}& =3i-800,\quad i=d/2+1,\ldots ,d, \\
m_{i}& =500+i,\quad i=d+1,\ldots ,2d,\\
m_{2d+1}& =1800.
\end{align*}
The internal damping is determined by \eqref{Dint1} with $\alpha =0.02$.
As shown on the Figure \ref{2doscillator}, we consider three
dampers. The first two dampers are grounded, while the third damper connects
two rows of masses. This means that external damping is determined by
\eqref{Dext} with
$$D_1= e_{i_1}e_{i_1}^T,\quad D_2=  e_{i_2}e_{i_2}^T\quad \mbox{and}\quad D_3=(e_{i_3}-e_{i_3+400})(e_{i_3}-e_{i_3+400})^T,$$
where $e_i$ corresponds to the $i$-th canonical vector.

In the definition of $D_{ext}$ we use indices since, in general, one  one needs to optimize damping positions, that is, we would like to optimize viscosities
$\rho_1, \rho_2$ and $\rho_3$ over different damping configurations.
Here we present results for only one configuration
$\begin{bmatrix} i_1 & i_2 & i_3\end{bmatrix}$.
We would like to damp 27 smallest eigenfrequencies of the undamped system, that is,
the matrix $G$ is defined by \eqref{G} with $s=27$.

In the homogeneous example $n=2001$,
$M=\mathop{\mathrm{diag}}(m_1,m_2,\ldots,m_n)$, where the first thousand masses are $m_i=1000$, the next thousand masses are $m_i=1500$, the last mass is $m_ {2001}=2000$, and
\begin{equation*}
K=\left[
\begin{array}{cccccc}
200 & -100 &  &  &  &  \\
-100 & 200 & -100 &  &  &  \\
& -100 & 200 &  &  &  \\
&  &  & \ddots &  &  \\
&  &  &  &  &  \\
&  &  &  & 300 & -150 \\
&  &  &  & -150 & 450
\end{array}
\right].
\end{equation*}
For this example we choose $s=20$ in \eqref{G}.

Our problems and solutions are described in Table \ref{tab:1}. The timings of the standard algorithm (Matlab)
and our new algorithm (Julia) are given in Tables \ref{tab:2} and \ref{tab:3}, respectively.

\begin{table}[H]
 \begin{tabular}{c|cccc}
  Problem & size & linearized size & $\begin{bmatrix} i_1 & i_2 & i_3\end{bmatrix}$ & optimal viscosities \\ \hline
  Small & 801 & 1602 & $\begin{bmatrix} 50 & 550 & 120\end{bmatrix}$ & $\begin{bmatrix} 561.4 & 651.8 & 310.6\end{bmatrix}$\\
  Large & 1601 & 3202 & $\begin{bmatrix} 50 & 950 & 220\end{bmatrix}$ & $\begin{bmatrix} 721.1 & 656.5 & 415.4\end{bmatrix}$ \\
  Homogeneous & 2001 & 4002 & $\begin{bmatrix} 850 & 1950 & 20\end{bmatrix}$ & $\begin{bmatrix} 620.0 & 1047.1 & 970.2\end{bmatrix}$
 \end{tabular}
 \caption{Size of test problem, size of linearized problem, indices defining configuration of dampers,
 and optimal viscosities. }
 \label{tab:1}
\end{table}

\begin{table}[H]
 \begin{tabular}{c|cc}
  Problem & lyap(SLICOT) & optimization \\ \hline
  Small & 1.8 & 162 (95 calls) \\
  Large & 11.2 & 1050 (97 calls) \\
  Homogeneous & 22.9 & 2608 (109 calls)
 \end{tabular}
 \caption{Run times in seconds for the standard $O(n^3)$ algorithm using Matlab with 12 cores. The individual
 problems are solved using
 function \texttt{lyap()} from SLICOT, and the unconstrained optimization is preformed using the function
 \texttt{fminsearchbnd()}.}
 \label{tab:2}
\end{table}

\begin{table}[H]
 \begin{tabular}{c|cccc}
  Problem & eigen(CSymDPR1) & traceX & optimization \\ \hline
  Small & 0.14 &  0.72 &  60 (79 calls)  \\
  Large & 0.48 &  2.4 &  212 (79 calls)  \\
  Homogeneous & 0.94 &  4.2 &  350 (71 calls)
 \end{tabular}
 \caption{Run times in seconds for the new $O(n^2)$ algorithm using Julia with 12 cores. The individual problems are
 solved with our eigensolver \texttt{eigen()} for CSYmDPR1 matrices followed by fast computation of trace with the
 function \texttt{traceX()}.
 The optimization is preformed using the function \texttt{ConjugateGradient()}
from the Julia package \texttt{Optim.jl} \cite{JuliaPack}.
Our Julia programs are available on GitHub \cite{Git}.}
 \label{tab:3}
\end{table}

We see that  the speedup of our algorithm over standard algorithm clearly grows with dimension (2.24, 4.03 and 4.8, computed by adjusting the number of trace computations).
Also, the computation times  in Table \ref{tab:3} for both, individual eigenvalue decomposition and trace computation, are clearly proportional to $n^2$.

To summarize, the proposed algorithm, based on CSymDPR1 matrix eigensolver and fast multiplication of Cauchy-like matrices, is simple, stable, and outperforms the standard counterpart, especially when the size of the problem is bigger. It is also easy to use in Julia's multithreading environment.

Some future work may include more detailed analysis of the the new eigensolver and development of the eigensolver for block complex symmetric diagonal-plus-low-rank matrices, which may treat all $k$ dampers simultaneously and, thus, be even faster.


\end{document}